\documentclass[MSNbibl,number,citesort,dvips]{arxbj}
\usepackage{upgreek}
\usepackage{graphicx}
\aid{0}
\volume{19}
\issue{1}
\pubyear{2013}
\firstpage{154}
\lastpage{171}
\doi{10.3150/11-BEJ397}

\makeatletter

\newtheorem{theo}{Theorem}

\newproclaim{defi}{Definition}
\newtheorem{cor}{Corollary}

\newremark{rem}{Remark}
\newcommand{\cal}{\mathcal}
\makeatother

\begin{document}
\begin{frontmatter}

\title{An extended family of circular distributions related to wrapped
Cauchy distributions via Brownian motion}
\runtitle{Circular distributions via Brownian motion}

\begin{aug}
\author[a]{\fnms{Shogo} \snm{Kato}\corref{}\thanksref{a}\ead[label=e1]{skato@ism.ac.jp}} \and
\author[b]{\fnms{M.C.} \snm{Jones}\thanksref{b}\ead[label=e2]{m.c.jones@open.ac.uk}}
\runauthor{S. Kato and M.C. Jones} 
\address[a]{Institute of Statistical
Mathematics, 10-3 Midori-Cho, Tachikawa, Tokyo 190-8562, Japan.\\\printead{e1}}
\address[b]{Department of Mathematics \& Statistics, The Open
University, Walton Hall, Milton Keynes MK7 6AA, UK. \printead{e2}}
\end{aug}

\received{\smonth{3} \syear{2011}}
\revised{\smonth{8} \syear{2011}}

\begin{abstract}
We introduce a four-parameter extended family of
distributions related to the wrapped Cauchy distribution on the
circle. The proposed family can be derived by altering the
settings of a problem in Brownian motion which generates the
wrapped Cauchy. The densities of this family have a closed form
and can be symmetric or asymmetric depending on the choice of the
parameters. Trigonometric moments are available, and they are shown to
have a simple form. Further tractable properties of the
model are obtained, many by utilizing the trigonometric moments.
Other topics related to the model, including alternative
derivations {and} M\"obius transformation,
are considered.
Discussion of the symmetric submodels is given. Finally,
generalization to a family of distributions on the sphere is
briefly made.
\end{abstract}


\end{frontmatter}

\section{Introduction}\label{sec1}

As a unimodal and symmetric model on the circle, the wrapped
Cauchy or circular Cauchy distribution has played an important role in
directional
statistics. It has density
%
\begin{equation} f(\theta;\mu,\rho) =
\frac{1}{2\uppi } \frac{1-\rho^2}{1+\rho^2 - 2\rho\cos(\theta-
\mu)}, \qquad-\uppi \leq\theta< \uppi , \label{wc}
\end{equation}
where $\mu \ (\in[-\uppi ,\uppi ))$ is a location parameter, and $\rho
\ ( \in[0,1))$ controls the concentration of the model.
This distribution has some desirable mathematical features, as
discussed in Kent and Tyler~\cite{ken} and McCullagh \cite
{mcc96}. We write
$\Theta\sim \operatorname{WC}(\mu,\rho)$ if a random variable
$\Theta$ has density (\ref{wc}).

In modeling symmetric circular data, the wrapped Cauchy
distribution could be one choice as well as other familiar
symmetric models such as the von Mises and wrapped normal
distributions; see also Jones and Pewsey~\cite{jon05}. In reality,
however, it is not so common to find such data in fields of
application. In many cases the data of interest are
asymmetrically distributed, and therefore probability
distributions having asymmetric densities are desired.

Construction of a tractable circular model with an asymmetric
shape has been a problem in statistics of circular data. To
tackle this problem, some asymmetric extensions of well-known
circular models have been proposed in the literature. Maksimov
\cite{mak}, Yfantis and Borgman~\cite{yfa} and Gatto and Jammalamadaka~\cite{gat}
discussed an extension of the von Mises distribution,
generated through maximization of Shannon's entropy with
restrictions on certain trigonometric moments.
Batschelet~\cite{bat} proposed a mathematical method of skewing
circular distributions that has seen renewed interest very
recently. Pewsey~\cite{pew} presented a four-parameter family
of distributions on the circle
by wrapping the stable distribution, which is asymmetric, onto
the circle. Recent work by Kato and Jones~\cite{kat10} proposed a
family of distributions arising from the M\"obius transformation
which includes the von Mises and wrapped Cauchy distributions.
Unlike familiar
symmetric distributions, it is often difficult to deal with skew
models in statistical analysis. This difficulty is partly due to
the lack of some mathematical properties that many of the
well-known symmetric models have. For example, existing
asymmetric models often have complex normalizing constants and
trigonometric moments, which could cause trouble in analysis.

In this paper, we provide a four-parameter extended family of
circular distributions, based on the wrapped Cauchy distribution,
by applying Brownian motion which, to our knowledge, has not been
used to propose a skew distribution on the circle. The four
parameters enable the model to describe not only symmetric shapes
but also asymmetric ones. An advantage of the proposed model is
its mathematical tractability. For instance, it has a simple
normalizing constant and trigonometric moments.
The
current proposal is more tractable than the family discussed by
Kato and Jones~\cite{kat10}, but is complementary to it in the sense
that the latter has other advantages (particularly some
associated with M\"obius transformation).

The subsequent sections are organized as follows. In Section~\ref{sec2} we
make the main proposal of this paper. The derivation of the
proposed model is given, and the probability density function and
probabilities of intervals under the density are discussed. Also,
special cases of the model are briefly considered. Section~\ref{sec3}
concerns the shapes of the density. Conditions for
symmetry and for unimodality are explored, and the interpretation
of the
parameters is discussed through pictures of the density. In
Section~\ref{sec4} we discuss the trigonometric moments and problems
related to them. It is shown that the trigonometric moments can
be expressed in a simple form. The mean direction,
mean resultant length and skewness of the model are also
considered.
Some other
topics concerning our model are provided in Section~\ref{sec5}. Apart from
the derivation given in Section~\ref{sec2}, there are other methods to
generate the family, which are discussed in that section. In
addition, we study the conformal transformation properties of the
distribution.
In Section~\ref{sec6} we investigate some
properties of the symmetric cases of the proposed model. These
submodels have some properties which the general family does not
have. Finally, in Section~\ref{sec7}, generalization to a family of
distributions on the sphere is made, and its properties are
briefly discussed.

Further material on the topic of this paper can be found in its technical
report
version, Kato and Jones~\cite{kat11}. This includes the path of
the expected exit points of a Brownian particle,
random variate generation using the Markov chain Monte Carlo method,
and proofs of Theorems~\ref{thm:density},~\ref{thm:symmetry}
and~\ref{thm:multi}.

\section{A family of distributions on the circle}\label{sec2}

\subsection{Definition}\label{sec2.1}
It is a well-known fact that the wrapped Cauchy distribution can be
generated as the distribution of the position of a Brownian
particle which
exits the unit disc in the two-dimensional plane; see, for example,
Section 1.10 in Durrett~\cite{dur}.
In addition to the wrapped Cauchy circular model, some distributions on
certain manifolds are derived from, or have relationships with, Brownian
motion. The Cauchy distribution on the real line is the distribution of
the position where a Brownian particle exits the upper half plane.
Considering the hitting time of the particle, one can obtain the inverse
Gaussian distribution on the positive half-line.
Kato~\cite{kat09} proposed a distribution for a pair of
unit vectors by
recording points where a Brownian particle hits circles with
different radii.

By applying Brownian motion, we provide a family of asymmetric
distributions on the circle, which includes the wrapped Cauchy
distribution as a
special case. The proposed model is defined as follows.
%
\begin{defi}
Let $\{B_t ; t \geq0\}$ be $\mathbb{R}^2$-valued Brownian motion
without drift
starting at $B_0= \rho_1 (\cos\mu_1, \sin\mu_1)'$, where $ 0 \leq
\rho_1 < 1$ and $-\uppi \leq\mu_1 < \uppi $. This Brownian particle will
eventually hit the unit circle. Let $\tau_1$ be the first time at
which the particle exits the circle, that is, $\tau_1 = \inf\{t ; \|
B_t\|
= 1\}$. After leaving the unit circle, the particle will hit a circle with
radius $\rho_2^{-1}\ (0 < \rho_2 <1)$ first at the time $\tau_2$, meaning
$ \tau_2=\inf\{t ; \|B_t\|= \rho_2^{-1}\}$. Then the proposed model is
defined by the conditional distribution of $B_{\tau_1}$ given $B_{\tau_2}=
\rho_2^{-1} (\cos\mu_2,\sin\mu_2)',$ where $-\uppi \leq\mu_2 < \uppi $.
\end{defi}

To put it another way, the proposed random vector, $B_{\tau_1}$ given
$B_{\tau_2}$, represents the position where a Brownian particle first
hits the unit circle, given the future point at which the particle
exits a circle with a larger radius.
From the next subsection, we investigate some properties of the
proposed model.

\subsection{Probability density function}\label{sec2.2}
One feature of the proposed model is that it has a closed form of
density with simple normalizing constant. It is given in the following
theorem.
See Appendix A of Kato and Jones~\cite{kat11} for the proof.

\begin{theo} \label{thm:density}
Let $\Theta$ be a $[-\uppi ,\uppi )$-valued random variable defined as
$B_{\tau_1} = (\cos\Theta, \sin\Theta)'$.
Then the conditional density of $\Theta$ given $B_{\tau_2}=\rho_2^{-1}
(\cos\mu_2,\sin\mu_2)'$ is given by
%
\begin{equation}
f(\theta) = C  [  \{ 1+\rho_1^2 -2\rho_1 \cos(\theta-\mu_1)
 \}  \{ 1+\rho_2^2 -2\rho_2 \cos(\theta-\mu_2)  \}
 ]^{-1} , \qquad-\uppi \leq\theta< \uppi , \label{density_theta}
\end{equation}
where $-\uppi \leq\mu_1,\mu_2 < \uppi , \ 0 \leq\rho_1,\rho_2 <1$, and the
normalizing constant $C$ is
\[
C= \frac{(1-\rho_1^2) (1-\rho_2^2) \{ 1+\rho_1^2 \rho_2^2 -2\rho_1 \rho_2
\cos(\mu_1-\mu_2) \}}{2\uppi (1-\rho_1^2 \rho_2^2)}.
\]
\end{theo}

For convenience, write $\Theta\sim \operatorname{EWC}(\mu_1,\mu_2,\rho_1,\rho_2)$ if a
random variable $\Theta$ has density (\ref{density_theta}) (``E''~for
``Extended'').
Note that the density does not involve any infinite sums or special
functions.

It is easy to see that distribution (\ref{density_theta}) reduces to the
wrapped Cauchy (\ref{wc}) when either $\rho_1$ or $\rho_2$ is equal to zero.

Quite often, it is advantageous to express the random variable
and
parameters of the proposed model in terms of complex numbers rather than
real numbers. Define a random vector by $Z=\mathrm{e}^{\mathrm{i} \Theta}$ where $\Theta
\sim \operatorname{EWC}(\mu_1,\mu_2,\rho_1,\rho_2)$.
Then $Z$ has density
%
\begin{equation}
f(z)= \frac{1}{2\uppi } \frac{|1-\phi_1 \overline{\phi_2}|^2}{1-|\phi_1
\overline{\phi_2}|^2} \frac{1-|\phi_1|^2} {|z-\phi_1|^2}
\frac{1-|\phi_2|^2}{|z-\phi_2|^2}, \qquad z \in\partial D, \label{density_z}
\end{equation}
with respect to arc length on the circle, where $\phi_1=\rho_1 \mathrm{e}^{\mathrm{i}
\mu_1},\ \phi_2 = \rho_2 \mathrm{e}^{\mathrm{i} \mu_2}$ and $\partial D= \{ z \in
\mathbb{C} ;  |z|=1 \}$.
It is clear that the parameters in this formulation, $\phi_1$
and
$\phi_2$, take values on the unit disc in the complex plane denoted by
$D=\{z \in\mathbb{C} ; |z|<1 \}$.
For brevity, we denote the distribution (\ref{density_z}) by $\operatorname{EC}^*
(\phi_1,\phi_2)$. Also, as in McCullagh~\cite{mcc96}, write $Z
\sim C^*(\phi_1)$
if $Z$ follows distribution (\ref{density_z}) with $\phi_2=0$.

We remark that the above change-of-variable and reparametrization do not
actually change the distribution. In later sections, we utilize either
of the representations (\ref{density_theta}) or (\ref{density_z}),
whichever is the more convenient.

\subsection{Probabilities}\label{sec2.3}
As well as the density of the proposed model, the probabilities of
intervals under the density can also be expressed without using infinite
sums or special functions.

\begin{theo} \label{thm:probabilities}
Let a random variable $\Theta$ follow the
$\operatorname{EWC}(\mu_1,\mu_2,\rho_1,\rho_2)$ distribution.
If $\rho_1 \neq\rho_2$ or $\mu_1 \neq\mu_2$, then the probability of
intervals under the density of $\Theta$ is given by
\begin{eqnarray*}
 &&P(a < \Theta\leq b)  \\
&& \quad = \int_{a}^{b} C  [  \{ 1+\rho_1^2 -2\rho_1 \cos(\theta-\mu
_1)  \}  \{ 1+\rho_2^2 -2\rho_2 \cos(\theta-\mu_2)  \}
 ]^{-1}\,\mathrm{d}\theta\\
&& \quad = \frac{C}{D}  \biggl\{ \rho_1 \rho_2 \sin(\mu_1-\mu_2)  \biggl[ \log
 \biggl\{ \frac{1+\rho_2^2-2 \rho_2 \cos(\theta-\mu_2)}{1+\rho_1^2-2
\rho_1 \cos(\theta-\mu_1)}  \biggr\}  \biggr]_a^b \\
&& \hphantom{\frac{C}{D}  \biggl\{}\qquad {} + \frac{ 2 \rho_1 \{ \rho_1 (1+\rho_2^2) - \rho_2 (1+\rho_1^2) \cos
(\mu_1-\mu_2) \} }{1-\rho_1^2}\\
&& \hphantom{{}+\frac{C}{D}  \biggl\{}\qquad {}\times  \biggl(  \biggl[ \arctan \biggl\{ \frac
{1+\rho_1}{1-\rho_1} \tan \biggl( \frac{\theta-\mu_1}{2}  \biggr)
\biggr\}  \biggr]_a^b + A_{\mu_1}  \biggr) \\
&& \hphantom{\frac{C}{D}  \biggl\{}\qquad {} + \frac{ 2 \rho_2 \{ \rho_2 (1+\rho_1^2) - \rho_1 (1+\rho_2^2) \cos
(\mu_1-\mu_2) \} }{1-\rho_2^2}\\
&& \hphantom{{}+\frac{C}{D}  \biggl\{}\qquad {}\times \biggl (  \biggl[ \arctan \biggl\{ \frac
{1+\rho_2}{1-\rho_2} \tan \biggl( \frac{\theta-\mu_2}{2}  \biggr)
\biggr\}  \biggr]_a^b + A_{\mu_2}  \biggr)  \biggr\},
\end{eqnarray*}
 where $C$ is defined as in Theorem~\ref{thm:density}, $-\uppi
\leq a<b<\uppi ,$
\[
D= (\rho_1^2 + \rho_2^2) (1+\rho_1^2 \rho_2^2) - 2 \rho_1 \rho_2
(1+\rho_1^2) (1+\rho_2^2) \cos(\mu_1-\mu_2) +4\rho_1^2\rho_2^2
\cos^2 (\mu_1-\mu_2)
\]
and
\[
A_{\mu} =
\cases{
{ \uppi ,} & \quad
$\displaystyle \tan \bigl\{ \tfrac12 (a-\mu)  \bigr\} > \tan \bigl\{ \tfrac12 (b-\mu)
 \bigr\} $,  \cr
{ 0,} & \quad  otherwise.
}
\]
If $\rho_1 = \rho_2$ and $\mu_1 = \mu_2$, then
\begin{eqnarray*}
P(a < \Theta\leq b)
&=& \int_a^b C  \{ 1+\rho_1^2 -2\rho_1 \cos (\theta-\mu_1)
 \}^{-2}\,\mathrm{d}\theta\\
&=& \frac{2 C}{(1-\rho_1^2)^2}  \biggl\{  \biggl[ \frac{\rho_1 \sin
(\theta-\mu_1)}{1+\rho_1^2-2\rho_1 \cos(\theta-\mu_1)}  \biggr]_a^b \\
&&\hphantom{\frac{2 C}{(1-\rho_1^2)^2}  \biggl\{}{} +
\frac{ 1+\rho_1^2}{1-\rho_1^2}  \biggl(  \biggl[ \arctan \biggl\{
\frac{1+\rho_1}{1-\rho_1} \tan
 \biggl( \frac{\theta-\mu_1}{2}  \biggr)  \biggr\}  \biggr]_a^b + A_{\mu_1}
 \biggr)  \biggr\}.
\end{eqnarray*}
\end{theo}

\begin{pf}
The result is straightforward from equations (2.553.3), (2.554.3) and
(2.559.2) of Gradshteyn and Ryzhik~\cite{gra}.
\end{pf}

\subsection{Special cases}\label{sec2.4}

The proposed family (\ref{density_theta}) contains some known
distributions as special cases.

\textit{Case 1:} The wrapped Cauchy distribution.
As mentioned in Section~\ref{sec2.2}, model (\ref{density_theta}) becomes the
wrapped Cauchy $\operatorname{WC}(\mu_1,\rho_1)$ if $\rho_2=0$.
Similarly, the model is $\operatorname{WC}(\mu_2,\rho_2)$ when $\rho_1=0$.

\textit{Case 2:} A special case of the Jones and Pewsey
\cite{jon05} family.
If $\rho_1=\rho_2  (\equiv\rho)$ and $\mu_1=\mu_2  (\equiv\mu)$, then
density (\ref{density_theta}) reduces to
%
\begin{equation}
f(\theta) = \frac{1}{2\uppi } \frac{1-\rho^2}{1+\rho^2}  \biggl(
\frac{1-\rho^2}{1+\rho^2-2\rho \cos(\theta-\mu)}  \biggr)^2, \qquad-\uppi
\leq\theta< \uppi . \label{d=2_c}
\end{equation}
This submodel corresponds to a special case of the family presented by
Jones and Pewsey~\cite{jon05}.
Their model has density
\[
f_{\mathrm{JP}}(\theta) = \frac{\{ \cosh(\kappa\psi) + \sinh(\kappa\psi) \cos
(\theta- \mu) \}^{1/\psi}}{ 2 \uppi  P_{1/\psi}(\cosh(\kappa\psi))},
\qquad-\uppi \leq\theta< \uppi ,
\]
where $\kappa\geq0,\ \psi\in\mathbb{R}, \ -\uppi \leq\mu<
\uppi $ and
$P_{1/\psi}(z)$ is the associated Legendre function of the
first kind of
degree $1/\psi$ and order 0 (Gradshteyn and Ryzhik~\cite{gra}, Sections
8.7, 8.8).
Our submodel is equivalent to their family with $\psi=-1/2$ and $\kappa=
2 \log((1+\rho)/(1-\rho))$.
In this case the associated Legendre function simplifies to
$P_{-2}(z)=P_1(z)=z$.
(See equations (8.2.1) and (8.4.3) of Abramowitz and Stegun~\cite{abr}.)
This fairly heavy-tailed circular distribution is shown by Jones and
Pewsey to
tend, suitably normalized, to the $t$ distribution on $\sqrt{3}$
degrees of
freedom as ($\rho\rightarrow1$ and hence) $\kappa\rightarrow\infty$!

\textit{Case 3:} One-point distribution.
As $\rho_1$ ($\rho_2$) tends to one, the model converges to a point
distribution with singularity at $\theta=\mu_1$ ($\theta=\mu_2$).
Normalizing by a scale factor of $(1-\rho)/\sqrt{\rho}$, the limiting
version of the wrapped Cauchy distribution as $\rho
\rightarrow1$ is the ordinary Cauchy distribution.
That argument can be
extended to show that the Cauchy limiting distribution also arises in
this case.

\textit{Case 4:} Two-point distribution.
Assume that $\rho_1=\rho_2$ $(\equiv\rho)$ and $\mu_1 \neq\mu_2$.
When $\rho$ goes to one, model (\ref{density_theta})
converges to the
distribution of a random variable which takes values on $\mu_1$ or
$\mu_2$ with probability 0.5.
If $\mu_1 < \mu_2$ and $\rho$ tends to one, $\{ \sqrt{\rho}/(1-\rho) \}
(\Theta-\mu_1) $ converges to $\frac12 C(0,1) + \frac12 I,$
where $C(0,1)$ is the standard Cauchy, and $I$ has the distribution
function $F(x) = \lim_{\lambda\rightarrow\infty} \int_{-\infty}^x 1/[
\uppi \{ 1+(t-\lambda)^2 \} ]\,\mathrm{d}t.$
Roughly speaking, this limiting distribution is a 50:50 mixture of the
standard Cauchy and a point distribution with singularity at $X=\infty$.
Similarly, one can show that if $\mu_1 > \mu_2$ and $\rho$ tends to
one, $\{ \sqrt{\rho}/(1-\rho) \} (\Theta-\mu_1)$ converges to a 50:50
mixture of the standard Cauchy and a point distribution with
singularity at $X=-\infty$.
The limiting distribution of $\{\sqrt{\rho} /(1-\rho)\} (\Theta-\mu_2)$
can be discussed in a similar manner.

\textit{Case 5:} The circular uniform distribution.
When $\rho_1=\rho_2=0$, the distribution reduces to the circular uniform
distribution.

\section{Shapes of probability density function}\label{sec3}

\subsection{Conditions for symmetry}\label{sec3.1}
As briefly stated in Section~\ref{sec1}, the proposed model can be symmetric or
asymmetric, depending on the values of the parameters.
The condition for symmetry is clearly written out in the following theorem.
%
\begin{theo} \label{thm:symmetry}
Density (\ref{density_theta}) is symmetric if and only if $\rho_1 =
\rho_2,\ \rho_j=0$ or $\mu_1=\mu_2+ (j-1) \uppi$ \mbox{$(j=1,2)$}.
\end{theo}

See Appendix B of Kato and Jones~\cite{kat11} for the proof.
The three-parameter symmetric submodels arising under the above
conditions
will be discussed in Section~\ref{sec6}.

\subsection{Conditions for unimodality}\label{sec3.2}
Turning to conditions for unimodality, it is possible to express these in
terms of an inequality.
The process to obtain the inequality is similar to that in Kato and Jones~\cite{kat10}, Section 2.5.
In the following discussion, take $\mu_2=0$ without loss of generality.
First we calculate the first derivative of density (\ref{density_theta})
with respect to $\theta$, which is given by
\[
\frac{\mathrm{d}}{\mathrm{d} \theta}f(\theta) \propto\frac{ [ \rho_1 \sin(\theta-\mu_1)
(1+\rho_2^2 -2\rho_2 \cos \theta) + \rho_2 \sin\theta\{ 1+\rho_1^2 -2
\rho_1 \cos(\theta-\mu_1) \} ] }{ \{1+\rho_1^2-2 \rho_1
\cos(\theta-\mu_1) \}^2 (1+\rho_2^2 -2 \rho_2 \cos\theta)^2 }.
\]
Then it follows that the extrema of the density are  obtainable as
solutions of the following equation:
%
\begin{equation}
a_0 +a_1 \cos\theta+ a_2 \sin\theta+ a_3 \cos\theta\sin\theta+
a_4 \cos^2 \theta=0, \label{extrema}
\end{equation}
where
\begin{eqnarray*}
a_0 &=& 2 \rho_1 \rho_2 \sin\mu_1, \qquad a_1 = \rho_1 (1+\rho_2^2) \sin
\mu_1, \qquad a_2 = -\rho_1 (1+\rho_2^2) \cos\mu_1 -\rho_2 (1+\rho
_1^2),\\
a_3 &=& 4 \rho_1 \rho_2 \cos\mu_1, \qquad a_4 = -4 \rho_1 \rho_2 \sin
\mu_1.
\end{eqnarray*}
Following Yfantis and Borgman~\cite{yfa} and Kato and Jones~\cite{kat10}, Section 2.5,
put $x = \tan ( \theta/2
 )$ so that $\cos\theta= (1-x^2)/(1+x^2)$, $\sin\theta= 2x
/(1+x^2)$. It follows that equation (\ref{extrema}) can be expressed as
a quartic equation in $x$ whose discriminant $D$ (Uspensky \cite
{usp}) can be
written down in
terms of $a_0,\ldots,a_4$ as in (6) of Kato and Jones~\cite{kat10}.
The quartic equation has four real roots or four complex
ones if $D>0$, and two real roots and two complex ones if $D<0$.
Therefore the distribution is bimodal when $D>0$ and unimodal when $D<0$.
Since $a_j \ (0 \leq j \leq4)$ are functions of $\mu_1,\rho_1$ and
$\rho_2$, we can write the conditions for unimodality as a function
of these three parameters. For general $\mu_2$, it is easy to see that the
discriminant is expressed in terms of $\mu_1-\mu_2$, $\rho_1$ and
$\rho_2$.

\begin{figure}[t]

\includegraphics{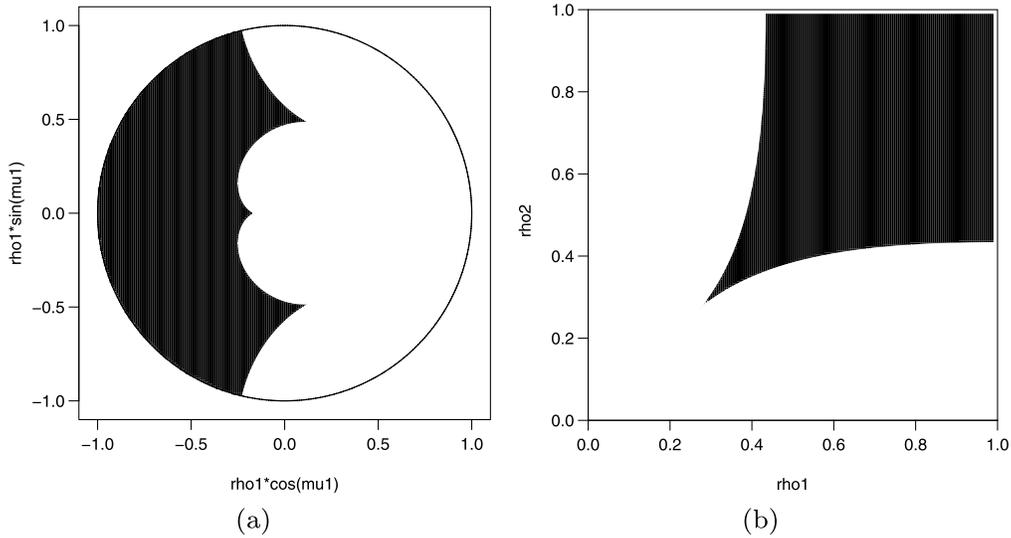}

\caption{Discriminant $D$
as functions of: (a) $(\rho_1 \cos\mu_1, \rho_1 \sin\mu_1)$
with $\rho_2=0.5$, and (b) $(\rho_1,\rho_2)$ with $\mu_1=2\uppi /3$.
Positive $D$, corresponding to bimodality, is displayed in black, and
negative $D$, corresponding to unimodality, is exhibited in
white.}
\label{fig1}\end{figure}

Figure~\ref{fig1} exhibits areas of
positivity (bimodality of density) and negativity (unimodality of
density) of the discriminant $D$ for two pairs of parameters.
Figure~\ref{fig1}(a) suggests that density (\ref{density_theta}) with
$\rho_1=0.5$ becomes unimodal for most cases when $-\uppi /2 \leq
\mu_1 \leq\uppi /2$. It seems that bimodality is most likely to
occur if $\uppi /2 \leq\mu_1 \leq3\uppi /2$. The other frame of
Figure~\ref{fig1} seems to show that distribution (\ref{density_theta})
with $\mu_2=2 \uppi /3$ is bimodal if $\rho_1$ or $\rho_2$ is
sufficiently large, and $\rho_2/\rho_1$ takes a value close to
$1$. In particular, when $\rho_1=\rho_2$, the range of parameters
corresponding to unimodality seems less wide than that when
$\rho_1=a \rho_2$ with $a\neq1$. As will be shown in Section~\ref{sec6.2}, for $\mu_1=2\uppi /3$ and $\rho_1=\rho_2$, the density becomes
unimodal if $\rho_1
\leq2 -\sqrt{3}\simeq0.268$ and bimodal otherwise.

\begin{figure}[t]

\includegraphics{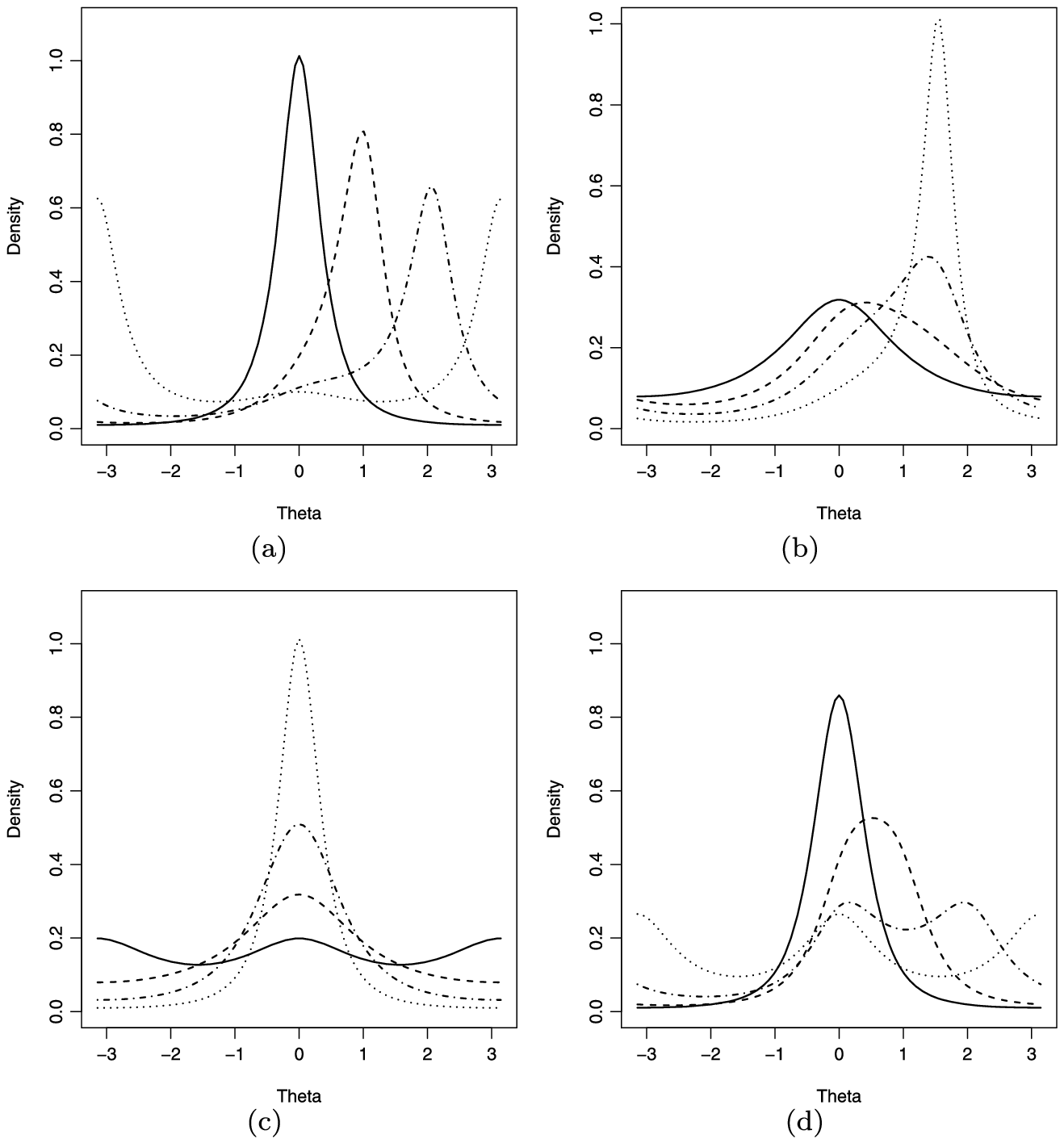}

\caption{Density (\protect\ref{density_theta}) for $\mu_2=0$ and
(a) $\rho_1=2/3,\ \rho_2=1/3$ and: $\mu_1=0$ (solid), $\uppi /3$ (dashed),
$2\uppi /3$ (dot-dashed) and $\uppi $ (dotted), (b) $\mu_1=\uppi /2,\ \rho_2=1/3$
and: $\rho_1=0$ (solid), $1/4$ (dashed), $1/2$ (dot-dashed) and $3/4$
(dotted), (c) $\rho_2=1/3$ and: $(\mu_1,\rho_1)=(\uppi ,1/3)$ (solid),
$(0,0)$ (dashed), $(0,1/3)$ (dot-dashed) and $(0,2/3)$ (dotted), and (d)
$\rho_1=\rho_2=0.5$ and: $\mu_1=0$ (solid), $\uppi /3$ (dashed), $2\uppi /3$
(dot-dashed) and $\uppi $ (dotted).}
\label{fig2}\end{figure}

\subsection{Pictures of density}\label{sec3.3}
We plot density (\ref{density_theta}) for
selected values of the parameters in Figure~\ref{fig2}.
Confirming the results in the previous two subsections, this figure shows
that the density can be symmetric or asymmetric, unimodal or bimodal,
depending on the choice of the parameters.
Figure~\ref{fig2}(a), displaying the density for fixed values of
$\mu_2,\rho_1$ and
$\rho_2$, suggests that the density is symmetric if $\mu_1=\mu_2$
or
$\mu_1=\mu_2+\uppi $ and asymmetric otherwise.
Figure~\ref{fig2}(b) exhibits how the shape of density (\ref{density_theta})
changes as $\rho_1$ increases.
It is clear from this frame that the greater the value of $\rho_1$, the
greater the concentration of the density.
The frame also implies that the absolute value of the circular skewness
of the density with fixed $\mu_1$ and $\mu_2$ takes a large value if
$\rho_1/\rho_2$ is close to 1.
Some symmetric cases of the density are shown in Figure~\ref{fig2}(c).
It appears to be that the density is unimodal when $\mu_1=\mu_2$, whereas
bimodality can occur when $\mu_2=\mu_1+\uppi $.
Figure~\ref{fig2}(d) is a case in which two parameters, $\rho_1$ and $\rho_2$, are
fixed to be equal (to 0.5), also implying that density
(\ref{density_theta}) is
symmetric. We will focus on this submodel in Section~\ref{sec6.2}.

\section{Trigonometric moments and related problems}\label{sec4}

\subsection{Trigonometric moments}\label{sec4.1}
It is often the case that asymmetric distributions on the circle have
complicated trigonometric moments or Fourier coefficients.
One feature of our asymmetric model, however, is the relative simplicity
of its trigonometric moments.
The expression for the moments is greatly simplified if the
variables and
parameters are represented in terms of complex numbers.

\begin{theo} \label{thm:tm}
Let a random variable $Z$ have the $\operatorname{EC}^*(\phi_1,\phi_2)$
distribution.
Then the $n$th trigonometric moment of $Z$ is given by
\[
E (Z^n) =
\cases{
\displaystyle  \frac{ (1-|\phi_2|^2) (1-\overline{\phi_1} \phi_2)
\phi_1^{n+1} - (1-|\phi_1|^2) (1-\phi_1 \overline{\phi_2})   \phi
_2^{n+1} }{(\phi_1 - \phi_2) ( 1 - | \phi_1 \overline{\phi_2} |^2 ) },
 & \quad
$\phi_1 \neq\phi_2$, \cr
\displaystyle \frac{1+n + (1-n) |\phi_1|^2}{1+ |\phi_1|^2} \phi_1^n ,
& \quad  $\phi_1 = \phi_2$.
}
\]
\end{theo}

\begin{pf}
The $n$th trigonometric moment can be expressed as
%
\begin{eqnarray}\label{tm_pr}
E(Z^n) &=& \frac{1}{2\uppi } \frac{|1-\overline{\phi_1} \phi_2|^2}{1-|\phi
_1 \overline{\phi_2}|^2} \int_{\partial D} z^n \Biggl ( \prod_{j=1}^2
\frac{1-|\phi_j|^2}{|z-\phi_j|^2}  \Biggr) |\mathrm{d}z| \nonumber
\\[-8pt]
\\[-8pt]
&=& \frac{1}{2\uppi  \mathrm{i}}   \frac{|1-\overline{\phi_1} \phi_2|^2}{1-|\phi
_1 \overline{\phi_2}|^2} \int_{\partial D} z^{n+1}  \Biggl\{ \prod
_{j=1}^2 \frac{1 - |\phi_j|^2}{(z-\phi_j)(1-\overline{\phi_j} z)}
\Biggr\}\,\mathrm{d}z .
\nonumber
\end{eqnarray}
If $\phi_1 \neq\phi_2$, the integrand in (\ref{tm_pr}) is holomorphic
on the unit disc, except at two poles of order 1, that is, $z=\phi_j,\ j=1,2$.
When $\phi_1 = \phi_2$, the integrand in (\ref{tm_pr}) has a single
pole of order 2 at $z=\phi_1$.
Hence, from the residue theorem (e.g., Rudin~\cite{rud}, Theorem 10.42), we
obtain Theorem~\ref{thm:tm}.
\end{pf}

\subsection{Mean direction and mean resultant length}\label{sec4.2}
Let $Z$ be distributed as $\operatorname{EC}^*(\phi_1,\phi_2)$.
From Theorem~\ref{thm:tm} it is easy to see that the first
trigonometric moment for $Z$ is simply expressed as
%
\begin{equation} \label{md}
E(Z) = \frac{(1-|\phi_2|^2) \phi_1 + (1-|\phi_1|^2) \phi_2}{ 1-|\phi_1
\phi_2|^2} \qquad\mbox{for any $\phi_1$ and $\phi_2$}.
\end{equation}
The following result provides the condition under which the first
trigonometric moment is equal to zero.
%
\begin{cor} \label{cor:tm=0}
The necessary and sufficient condition for $E(Z)=0$ is given by $\phi_1 =
-\phi_2$.
\end{cor}
\begin{pf}
It is straightforward to show the sufficient condition for $E(Z)=0$.
The necessary condition is proved as follows.
Let $E(Z)=0$ and $\phi_1 \neq- \phi_2$.
Then, from Theorem~\ref{thm:tm}, the following equation holds between the
two parameters:
$
\phi_1 \phi_2 = ({\phi_1 + \phi_2})/{ \overline{(\phi_1+\phi_2)}}.
$
Taking the absolute values of both sides of the above equation, we have
$
|\phi_1 \phi_2| =1.
$
Since $|\phi_j|<1 \ (j=1,2)$, there does not exist $\phi_1$ and $\phi
_2$ such that $|\phi_1 \phi_2|=1$.
Therefore, if $\phi_1 \neq-\phi_2$, then $E(Z) \neq0$.
\end{pf}

\begin{figure}

\includegraphics{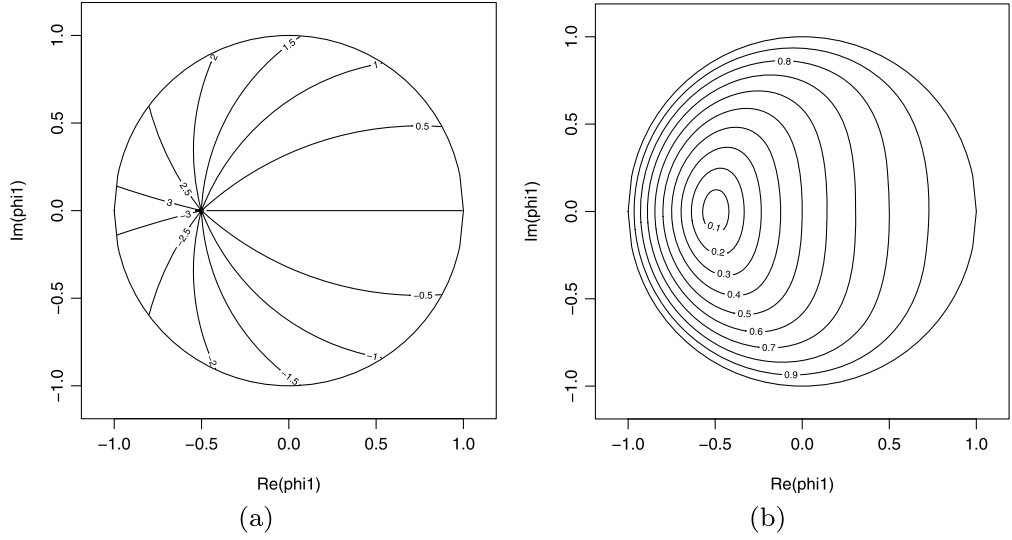}

\caption{Plot of (a) mean direction and (b) mean resultant
length of a random variable having the distribution
$\operatorname{EC}^*(\phi_1,0.5)$, as
a function of $\phi_1$.}
\label{fig3}\end{figure}

The mean direction is a measure of location which is defined as $\arg
\{E(Z) \}$ if $E(Z) \neq0$ and is undefined if $E(Z)=0$.
The mean resultant length, which is a measure of concentration, is
defined by $|E(Z)|$. Figure~\ref{fig3} shows the mean direction and mean resultant
length of a random variable having density (\ref{density_z}) with fixed
$\phi_2=0.5$.
Figure~\ref{fig3}(a) suggests that the mean direction monotonically increases as
$\arg(\phi_1)$ increases from $-\uppi $ to $\uppi $.
As seen in Corollary~\ref{cor:tm=0}, Figure~\ref{fig3}(b) confirms that the mean
resultant length is zero if $\phi_1=-0.5$.
It seems from this frame that unimodality holds for the mean
resultant length as a function of $\phi_1$.

\subsection{Skewness}\label{sec4.3}

A measure of skewness for circular distributions (Mardia~\cite{mar}), $s$, is
defined by
%
\begin{equation}
s = {E[\operatorname{Im}\{ (Z \mathrm{e}^{-\mathrm{i} \zeta})^2 \}]} / {(1-\delta)^{3/2}},
\label{skewness}
\end{equation}
where $\zeta$ and $\delta$ are the mean direction and mean resultant length
for a $\partial D$-valued random variable $Z$, respectively.
For our model, it is possible to express the skewness for the proposed
model in a fairly simple form as follows.
%
\begin{cor} \label{cor:skew}
Let $Z$ be distributed as $\operatorname{EC}^*(\phi_1,\phi_2)$.
Then the skewness for the distribution of $Z$ is given by
%
\begin{equation}
s = \frac{|1-\phi_1 \overline{\phi_2}|^2}{ ( 1-|\phi_1 \overline{\phi
_2}|^2)^3} \rho^{-2} (1-\rho)^{-3/2}   \operatorname{Im}  ( \phi_1
\overline{\phi_2}  ) (|\phi_1|^2-|\phi_2|^2) (1-|\phi_1|^2)
(1-|\phi_2|^2), \label{skew}
\end{equation}
where $\rho=  | (1-|\phi_2|^2) \phi_1 + (1-|\phi_1|^2) \phi_2
 | / (1-|\phi_1 \overline{\phi_2}|^2).$
\end{cor}
\begin{pf}
Note that the numerator of (\ref{skewness}) can be expressed as
\[
E[\operatorname{Im}\{ (Z \mathrm{e}^{-\mathrm{i} \mu})^2 \}] = \operatorname{Im}  [ E (Z^2)   \{
\overline{ E(Z) }  / |E(Z)|  \}^2  ].
\]
Then, by using Theorem~\ref{thm:tm} and equation (\ref{md}), the
skewness (\ref{skew}) can be obtained after lengthy, but
straightforward, calculations.
\end{pf}

Now we write $s=s(\phi_1,\phi_2)$ in order to view the skewness
(\ref{skew}) as a function of $\phi_1$ and $\phi_2$.
It follows from Corollary~\ref{cor:skew} and Theorem~\ref{thm:symmetry}
that the following properties hold for $s$.
%
\begin{theo} \label{thm:s}
The following properties hold for the skewness (\ref{skew}):
\begin{enumerate}[5.]
\item$s(\phi_1,\phi_2) = 0  \Longleftrightarrow$   density
(\ref{density_z}) is symmetric.
\item$s(\phi_1,\phi_2) > 0 \ (<0)  \Longleftrightarrow \operatorname
{Im}(\phi_1 \overline{\phi_2}) (|\phi_1| - |\phi_2|) > 0 \ (<0).$
\item$s(\phi_1,\phi_2) = -s (\phi_1, \overline{\phi_2} \phi_1^2 / |\phi
_1|^2) = - s (\overline{\phi_1} \phi_2^2 / |\phi_2|^2,\phi_2),
\phi_1,\phi_2 \neq0$.
\item$s(\phi_2,\phi_1)=s(\phi_1,\phi_2)$.
\item$s(\overline{\phi_1},\overline{\phi_2})=-s(\phi_1,\phi_2).$
\item$s(\alpha\phi_1,\alpha\phi_2) = s(\phi_1,\phi_2),  \alpha
\in\partial D.$
\item Assume $\rho_1,\rho_2>0$ and $\sin(\mu_1-\mu_2)>0\ (<0)$.
Then $\lim_{\rho_1 \rightarrow1} s(\rho_1 \mathrm{e}^{\mathrm{i} \mu_1},\rho_2 \mathrm{e}^{\mathrm{i} \mu
_2}) = \infty\ (-\infty)$.
\end{enumerate}
\end{theo}

\begin{figure}[b]

\includegraphics{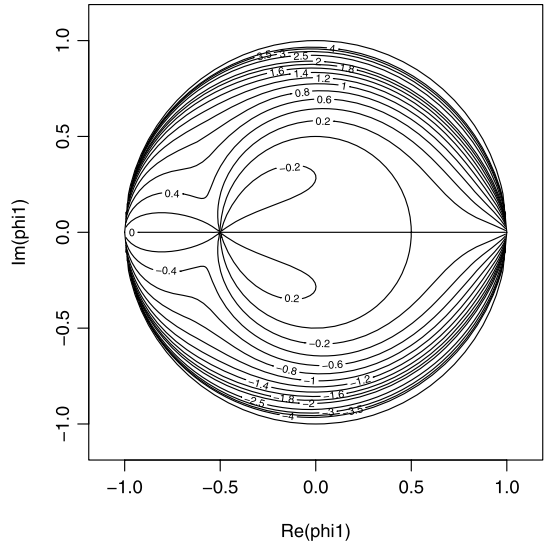}

\caption{Skewness (\protect\ref{skew}) of a random variable
having the distribution $\operatorname{EC}^*(\phi_1,0.5)$, as a function of
$\phi_1$.}
\label{fig4}\end{figure}

Figure~\ref{fig4} plots this skewness when $\phi_2=0.5$.
As the first property of Theorem~\ref{thm:s} shows, the skewness is
equal to zero if and
only if $\operatorname{Im}(\phi_1 \overline{\phi_2})=0,\ |\phi_1|=|\phi_2|$ or
$|\phi_j|=0$ $(j=1,2)$.
This figure also confirms the second property of Theorem~\ref{thm:s}
that the model is positively (negatively)
skewed if $|\phi_1|>|\phi_2|$ and $\operatorname{Im}(\phi_1 \overline{\phi_2})>0$
($\operatorname{Im}(\phi_1 \overline{\phi_2})<0$), or $|\phi_1|<|\phi_2|$ and
$\operatorname{Im}(\phi_1 \overline{\phi_2})<0$ ($\operatorname{Im}(\phi_1
\overline{\phi_2})>0$).
It can also be seen that the fifth and seventh properties of
Theorem~\ref{thm:s} hold in this figure.

\section{Some other topics}\label{sec5}

\subsection{Alternative derivations}\label{sec5.1}

As discussed in Section~\ref{sec2.1}, the proposed model can be derived by considering
a problem in Brownian motion. Apart from the derivation given there,
here is
another method to generate our family via this stochastic process, kindly
suggested by the referee.

Remember that, given a $\mathbb{C}$-valued Brownian path without drift
starting from $\phi$ $(\in D)$, the point of first exit from the unit
disc $D$
has a wrapped Cauchy distribution with density proportional to
$1/|z-\phi|^2$.
(Although we assume $\phi\in D$ by convention, the Brownian motion starting
at $1/\overline{\phi}$ has the same hitting distribution from the
exterior because of the
conformal invariance of Brownian motion.) From this it follows that another
derivation of the proposed model (\ref{density_z}) is given in the following
theorem.
%
\begin{theo} \label{thm:alt_br} Let $Z_1$ and $Z_2$ be the points of
first exit from the unit disc for two independent $\mathbb{C}$-valued Brownian
motions without drift starting from $\phi_1$ and $\phi_2$,
respectively, where
$\phi_1,\phi_2 \in D$. Then the conditional distribution of $Z_1$
$(\equiv Z)$
given $Z_1=Z_2$ is given by density (\ref{density_z}).
\end{theo}

%
\begin{rem}  The referee embedded the above in a generalized formulation
leading to densities of the form
%
\begin{equation} f (z) \propto\prod_{j=1}^d
\frac{1-|\phi_j|^2}{|z - \phi_j|^2}, \qquad z \in\partial D;
\phi_1,\ldots,\phi_d \in D, \label{extension}
\end{equation}
many of whose
properties are also readily obtainable using the residue theorem.
\end{rem}

In addition to the derivations using Brownian motion, there are other
methods to
generate our family.
The following result shows that the proposed model appears in a Bayesian
context.

\begin{theo} \label{thm:bayesian}
Let $\Theta| \nu$ have the wrapped Cauchy distribution $\operatorname{WC}(\nu,\rho_1)$
with known $\rho_1$.
Assume that the prior distribution of $\nu$ is distributed as the wrapped
Cauchy $\operatorname{WC}(\mu_2,\rho_2)$.
Then the posterior distribution of $\nu$ given $\Theta= \mu_1$ has
density (\ref{density_theta}).
\end{theo}

The above derivation enables us to generate random samples from our model
using the Markov chain Monte Carlo method;
see Kato and Jones~\cite{kat11}, Section 5.3, for details.

\begin{rem}
In Theorem~\ref{thm:bayesian}, the marginal
distribution of $\Theta$ is given by the wrapped Cauchy
distribution,
$\operatorname{WC}(\mu_1+\mu_2,\rho_1 \rho_2)$.
\end{rem}

Before we discuss a third method to derive our model, here we briefly
recall
a known property of the wrapped Cauchy distribution or the Poisson kernel
as follows.
%
\begin{theo}[(Rudin~\cite{rud}, Theorem 11.9)] \label{thm:rudin}
Assume $\overline{D}$ is the closed unit disc in the complex plane.
Suppose a function $u$ is continuous on $\overline{D}$ and harmonic in
$D$, and suppose that $|u(z)|<\infty$ for any $z \in\overline{D}$.
Then
\[
\int_{\partial D} u(z) \frac{1}{2 \uppi } \frac{1-|\phi_1|^2}{|z-\phi
_1|^2}\,\mathrm{d}z = u(\phi_1), \qquad\phi_1 \in D.
\]
\end{theo}

Note that the integrand in the left-hand side of the above equation is
actually the product of $u(z)$ and the density of the wrapped Cauchy
$C^*(\phi_1)$. From the fact that
\[
g(z) = \frac{|\phi_2^{-1}|^{2} - |z|^2}{ |z- \overline{\phi_2}^{-1}
 |^2}, \qquad\phi_2 \in D,
\]
is a continuous real function on $\overline{D}$ and harmonic in $D$, it
follows that
\[
f(z)= \frac{g(z)}{2 \uppi  g(\phi_1)} \frac{1-|\phi_1|^2}{|z-\phi_1|^2},
\qquad z \in\partial D,
\]
satisfies the definition of a density.
Clearly, $f(z)$ is equal to density (\ref{density_z}).

The following theorem implies that our model appears as a special case of
a further general family of distributions.
The proof is straightforward and omitted.
%
\begin{theo}
Let $f$ and $g$ be probability density functions on the circle.
Assume that $h$ is the convolution of $f$ and $g$, namely, $h(\tau
)=(f*g)(\tau)=\int_{-\uppi }^{\uppi } f(\theta) g(\tau-\theta)\,\mathrm{d}\theta$.
Then a function $k(\theta)=f(\theta) g(\tau-\theta) /h(\tau)$ is a
probability density function on the circle.
\end{theo}

Density (\ref{density_theta}) can be derived on setting $\tau=\mu_2$ and
$f$ and $g$ as the densities of $\operatorname{WC}(\mu_1,\rho_1)$ and $\operatorname{WC}(0,\rho_2)$,
respectively.
In this case $h(\tau)$ also has the form of the wrapped Cauchy density
$\operatorname{WC}(\mu_1,\rho_1 \rho_2)$ since the wrapped Cauchy distribution has the
additive property.

\subsection{M\"obius transformation}\label{sec5.2}
The M\"obius transformation is well known in complex analysis as a
conformal mapping which projects the unit disc onto itself.
It is defined as
\[
{\cal{M}}(w) = \alpha  \frac{w +\beta}{\overline{\beta} w +1 }, \qquad w
\in\overline{D};   \alpha\in\partial D,  \beta\in D.
\]
Although this projection is usually defined on the interior of the unit
disc, here we extend the domain of the mapping so that the boundary of
the unit disc, that is, the unit circle, is included.
It is easy to see that the unit circle is mapped onto itself via the
transformation, namely, ${\cal M}(\partial D) = \partial D$.
In directional statistics the transformation appears in some papers such
as McCullagh~\cite{mcc96}, Jones~\cite{jon04} and
Kato and Jones~\cite{kat10}.

Here we study the conformal invariance properties of our family.
The distribution (\ref{density_z}) has a density that is a relative
invariant of weight 2 in the sense that
\[
f(z;\phi_1,\phi_2) = f \{ {\cal M}(z) ; {\cal M}(\phi_1), {\cal M}(\phi
_2) \}   |{\cal M}'(z)|^2.
\]
Similarly, the extended family (\ref{extension}) has a density that is
relatively invariant with weight $k$, in the sense that $ f(z;\phi
_1,\ldots,\phi_k) = f \{ {\cal M}(z) ; {\cal M}(\phi_1), \ldots, {\cal
M}(\phi_k) \}   |{\cal M}'(z)|^k. $
The relative invariance with weight 1 of the wrapped Cauchy
distribution can be obtained by putting $k=1$.

The proposed family (\ref{density_z}) is not closed under the M\"obius
transformation except for the wrapped Cauchy special cases.
To say the same thing in a different way, if $Z$ follows the
distribution (\ref{density_z}), then the density for ${\cal M}(Z)$ is
not of the form (\ref{density_z}) except for $\phi_1=0$ or $\phi_2=0$.
This fact can be understood clearly in the following context;
as seen in Theorem~\ref{thm:alt_br}, our model (\ref{density_z}) can be
derived as the conditional distribution of $Z_1$ given $Z_1=Z_2$, where
$Z_1$ and $Z_2$ are independent wrapped Cauchy variables.
Invertibility of the transformations ${\cal M}$ implies that $Z_1=Z_2$
if and only if ${\cal M}(Z_1) = {\cal M}(Z_2)$.
However, as is known as the Borel paradox (e.g., Pollard~\cite{pol}, Section 5.5), the conditional distribution of
$Z_1$ given $Z_1=Z_2$ is generally not the same as the conditional
distribution of $Z_1$ given ${\cal M}(Z_1) = {\cal M}(Z_2)$.
Consequently, the M\"obius transformation of the conditional
distribution of $Z_1$ given $Z_1=Z_2$, which is the M\"obius
transformation of model (\ref{density_z}), is not the same as the
conditional distribution of ${\cal M}(Z_1)$ given ${\cal M}(Z_1)={\cal
M}(Z_2)$, which is model (\ref{density_z}).
In a similar manner one can see that the extended model (\ref{extension})
with $k \geq2$ is not closed under the M\"obius transformation.

\section{Symmetric cases}\label{sec6}

So far, we have mostly considered the full family of
distributions with densities
(\ref{density_theta}) and (\ref{density_z}).
In this section, we focus on the symmetric special cases of the proposed
model. Some properties, which the general family does not have, hold
for the
symmetric cases.
In Section~\ref{sec6.1}, model (\ref{density_theta}) with
$\mu_1=\mu_2$ or
$\mu_2+\uppi $ is discussed.
In Section~\ref{sec6.2}, we briefly consider another symmetric model,
namely,
model (\ref{density_theta}) with $\rho_1=\rho_2$.

\subsection{\texorpdfstring{Symmetric case I: $\mu_1=\mu_2$ or $\mu_1=\mu_2+\uppi $}
{Symmetric case I: mu_1=mu_2 or mu_1=mu_2+pi}}\label{sec6.1}
Model (\ref{density_theta}) with $\mu_1=\mu_2$ or
$\mu_1=\mu_2+\uppi $
is essentially the same as the distribution with density
%
\begin{eqnarray}\label{symmetric1}
&&f(\theta) = \frac{1-\rho_1 \rho_2 }{2 \uppi (1+\rho_1 \rho_2)}
\frac{ 1-\rho_1^2 }{ 1+\rho_1^2 -2 \rho_1 \cos(\theta- \mu) }
\frac{1-\rho_2^2 }{ 1+\rho_2^2 -2 \rho_2 \cos(\theta- \mu) },\nonumber
\\[-8pt]
\\[-8pt]
&& \quad
-\uppi \leq\theta< \uppi ,
\nonumber
\end{eqnarray}
where $-1 < \rho_1 <1, \ 0 \leq\rho_2<1,$ and $-\uppi \leq\mu< \uppi $.
Note the extension in the range of $\rho$:
this density corresponds to (\ref{density_theta}) with $\mu_1=\mu_2$ when
$\rho_1 \geq0$ in (\ref{symmetric1}), while the density is equivalent
to (\ref{density_theta}) with $\mu_1=\mu_2+\uppi $ when $\rho_1 <0$ in
(\ref{symmetric1}).
Clearly, the density is symmetric about $\theta=\mu$ and $\theta= \mu
+\uppi $.

One might, however, wish to restrict interest to the case where
$0\leq\rho_1<1$ because then density (\ref{symmetric1}) is unimodal. To
see this, note that in this case $a_0 = a_1 = a_4=0$ in (\ref{extrema})
so that stationary points of the density occur at $\sin\theta= 0$ and,
potentially, also at $\cos\theta= -a_2/a_3 =
(\rho_1+\rho_2)(1+\rho_1\rho_2)/4\rho_1\rho_2$. However, the
latter expression can easily be proved to not be less than unity, showing
that (\ref{symmetric1}) is unimodal.
On the other hand, the symmetric model (\ref{symmetric1}) with negative
$\rho_1$ has a relationship with mixtures of two wrapped Cauchy
distributions.
The density for this submodel can be expressed as
\begin{eqnarray*}
f(\theta) &=& p \frac{1}{2\uppi } \frac{1-\rho_1^2}{1 + \rho_1^2 -2\rho_1
\cos(\theta-\mu) } + (1-p) \frac{1}{2\uppi } \frac{1-\rho_2^2}{1 + \rho
_2^2 -2\rho_2 \cos(\theta -\mu) },
\end{eqnarray*}
where $p= -\rho_1 (1-\rho_2^2)/ \{(1+\rho_1 \rho_2)(\rho_2 - \rho_1)\}$.
Inter alia, this representation allows random numbers from density
(\ref{symmetric1}) with $\rho_1 \leq0$ to be more
easily generated than they were for the general case in Kato and Jones~\cite{kat11},
Section 5.3.

\subsection{\texorpdfstring{Symmetric case II: $\rho_1=\rho_2$}{Symmetric case II: rho_1=rho_2}}\label{sec6.2} In this section
we discuss the submodel of density (\ref{density_theta})
with
$\rho_1=\rho_2=\rho \ (\in[0,1))$. The density for this
submodel reduces to
\begin{eqnarray} f(\theta) &=&
\frac{(1-\rho^2) \{ 1+\rho^4-2\rho^2 \cos(\mu_1-\mu_2) \}}{2\uppi
(1+\rho^2)} \nonumber\\
&&{} \times [  \{ 1+\rho^2 -2
\rho\cos(\theta- \mu_1)  \}  \{ 1+\rho^2 -2 \rho\cos
(\theta- \mu_2)  \}  ]^{-1}, \qquad-\uppi \leq\theta<
\uppi . \nonumber
\end{eqnarray}
The model is symmetric about
$\theta=\overline{\mu}$ and $\overline{\mu}+\uppi $, where
$\overline{\mu} = \arg( \sum_{j=1}^2 \cos\mu_j + \mathrm{i} \sum_{j=1}^2
\sin\mu_j )$. As displayed in Figure~\ref{fig2}(d), this submodel can be
unimodal or bimodal depending on the choice of the parameters.
The condition for unimodality can be simplified for this submodel
as follows:
\[
-2 \arccos \biggl( \frac{2 \rho}{1+\rho^2} \biggr )
\leq\mu_1 -\mu_2 \leq2 \arccos \biggl( \frac{2 \rho}{1+\rho^2}
\biggr ).
\]
This result was used in Section~\ref{sec3.2} to find the
region of Figure~\ref{fig1}(b)
in which the discriminant takes negative values, corresponding to
unimodality, when $\rho_1=\rho_2$ and $\mu_1-\mu_2=2\uppi /3$.
\section{A generalization on the sphere}\label{sec7}

As described in Section~\ref{sec2.1}, the proposed model
(\ref{density_theta}) can be derived using Brownian motion.
By adopting a multi-dimensional Brownian motion instead of the
two-dimensional one, we can extend model
(\ref{density_theta}) to a distribution on the unit sphere.
The generalized model is defined as follows.

\begin{defi}
Let $\{B_t ; t \geq0\}$ be $\mathbb{R}^d$-valued Brownian motion
starting at $B_0= \rho_1 \eta_1$, where \mbox{$d \geq2$}, $0 \leq\rho_1 < 1
,\ \eta_1 \in S^d $ and $S^d = \{ x \in\mathbb{R}^d ; \|x\|=1 \}$.
Assume $\tau_1 = \inf\{t ; \|B_t\|= 1\}$ and $ \tau_2=\inf\{t ; \|
B_t\|= \rho_2^{-1}\}$, where $0 < \rho_2 <1$.
Then the proposed model is defined by the conditional distribution of
$B_{\tau_1}$ given $B_{\tau_2}= \rho_2^{-1} \eta_2,$ where $\eta_2 \in S^d$.
\end{defi}

For simplicity, write $X=B_{\tau_1}$.
The probability density function of this extended model is available,
and it is given in the following theorem.
%
\begin{theo} \label{thm:multi}
The conditional distribution of $X$ given $B_{\tau_2}=\rho_2^{-1} \eta
_2$ is of the form
%
\begin{equation}
f(x) = \frac{1}{A_{d-1}} \frac{(1 +\rho_1^2 \rho_2^2 -2 \rho_1 \rho_2
\eta_1' \eta_2)^{d/2}}{1-\rho_1^2 \rho_2^2} \frac{1-\rho_1^2}{\|x-\rho
_1 \eta_1\|^d} \frac{1-\rho_2^2}{\|x-\rho_2 \eta_2\|^d}, \qquad x \in
S^d, \label{sphere}
\end{equation}
where $A_{d-1}$ is the surface area of $S^d$, that is, $A_{d-1}=2 \uppi
^{d/2} / \Gamma(d/2) $.
\end{theo}

See Appendix C of Kato and Jones~\cite{kat11} for the proof.
Note that density (\ref{sphere}) reduces to the circular case (\ref
{density_theta}) if $d=2,\ x=(\cos\theta, \sin\theta)',\ \eta_1=(\cos
\mu_1,\sin\mu_1)'$ and $\eta_2 = (\cos\mu_2,\sin\mu_2)'$.
It might be appealing that density (\ref{sphere}), which is not
rotationally symmetric in general, can be expressed in a relatively
simple form.

When $\rho_2=0$, the distribution becomes the so-called ``exit''
distribution on the sphere, whose density is given by
%
\begin{equation}
f(x) = \frac{1}{A_{d-1}} \frac{1-\rho_1^2}{\|x - \rho_1 \eta_1\|^d},
\qquad x \in S^d. \label{exit}
\end{equation}
We write $X \sim\operatorname{Exit}_d (\rho_1 \eta_1)$ if a random variable
$X$ has density (\ref{exit}).
This model is rotationally symmetric about $x=\eta_1$.
See, for example, Durrett~\cite{dur}, Section 1.10, for details about the distribution.
In particular, when $d=2$, model (\ref{exit}) becomes the wrapped
Cauchy distribution.
It is noted that this distribution is a submodel of Jones and Pewsey's
\cite{jon05} family of distributions on the sphere with density
%
\begin{equation}
f(x) = \frac{1}{A_{d-1}} \frac{|\sinh(\kappa\psi)|^{d/2-1}}{2^{2/d-1}
\Gamma(d/2)} \frac{ \{\cosh(\kappa\psi) + \sinh(\kappa\psi) x' \mu
\}^{1/\psi}}{P_{1/\psi+d/2-1}^{1-d/2} \{ \cosh(\kappa\psi)\}}, \qquad
x \in S^d, \label{jp_sphere}
\end{equation}
where $P_{\alpha}^{\beta} (z)$ is the associated Legendre function of
the first kind of degree $\alpha$ and order $\beta$ (Gradshteyn and Ryzhik~\cite{gra}, Sections 8.7, 8.8).
It follows from equation (8.711.1) of Gradshteyn and Ryzhik~\cite{gra} and
equation (2) of McCullagh~\cite{mcc89} that density
(\ref{jp_sphere}) reduces to density (\ref{exit}) if $\kappa= d
\log\{ (1+\rho_1)/(1-\rho_1) \}/2, \ \mu=\eta_1$ and $\psi=-2/d$.
Also, the extended model (\ref{sphere}) includes the model with density
%
\begin{equation}
f(x) = \frac{1}{A_{d-1}} \frac{(1-\rho_1^2)^{d+1}}{1+\rho_1^2} \frac
{1}{\| x - \rho_1 \eta_1 \|^{2d}}, \qquad x \in S^d, \label{d=2_s}
\end{equation}
when $\rho_1=\rho_2,\ \eta_1=\eta_2$.
It can be seen that this model is another submodel of Jones
and Pewsey's~\cite{jon05} family (\ref{jp_sphere}) by putting
$\psi=-1/d$ and $\kappa=d \log\{(1+\rho_1)/(1-\rho_1)\}$.
Also, notice that, if $d=2$, density (\ref{d=2_s})
corresponds to Case 2, (\ref{d=2_c}), of circular submodels in
Section~\ref{sec2.4}.
In addition to these two submodels, the multivariate model
(\ref{sphere}) contains the uniform distribution ($\rho_1=\rho_2=0$),
one-point distribution ($\rho_1=0$ and $\rho_2 \rightarrow1$, or $\rho
_2 =0$ and $\rho_1 \rightarrow1$), and two-point distribution ($\rho
_1=\rho_2,\ \mu_1 \neq\mu_2$ and $\rho_1 \rightarrow1$).

In a similar manner to Theorem~\ref{thm:bayesian}, it is easy to
prove that model (\ref{sphere}) can be derived from Bayesian
analysis of the exit distribution.

\begin{theo} \label{thm:bayesian_sphere}
Let $X | \xi$ be distributed as the exit distribution $\operatorname{Exit}(\rho
_1 \xi)$ with known $\rho_1$.
Suppose that the prior distribution of $\xi$ is
$\operatorname{Exit}(\rho_2 \eta_2)$.
Then the posterior distribution of $\xi$ given $X=\eta_1$ is
given by density (\ref{sphere}).
\end{theo}

From this result and the random number generator of the exit
distribution given by Kato~\cite{kat09}, Section 3.3, it is possible
to generate random samples from (\ref{sphere}) through the
Markov chain Monte Carlo method in a similar manner to Kato and Jones~\cite{kat11},
Section~5.3.

\section*{Acknowledgements}
The authors are grateful to the Associate Editor and a referee for suggesting
a simpler proof of Theorem~\ref{thm:tm}, the alternative derivation of the
proposed model via Brownian motion given in Theorem~\ref{thm:alt_br}
and the
conformal invariance properties discussed in Section~\ref{sec5.2}.
Financial support for this research was received by Kato in the form of
``Grant-in-Aid for Young Scientists (B)'' (22740076) from Japan Society
for the Promotion of Science.


\printhistory

\end{document}